# The *n*-dimensional Cube---A New Way to Prove the Fermat's Last Theorem

**Kaida Shi**

State Key Laboratory of CAD&CG, Zhejiang University,
Hangzhou 310027, Zhejiang Province, China
Department of Mathematics, Zhejiang Ocean University,
Zhoushan 316004, Zhejiang Province, China

**Abstract** By researching the *n*-dimensional cube, the author found a new way, and proved that the indefinite equation $x^n + y^n = z^n \,(n \geq 3)$ has no any solutions of positive integers.

The famous hypothesis "when *n*≥3, the indeterminate equation $x^n + y^n = z^n$ has no any solutions of positive integers", raised by French mathematician Pierre de Fermat in 1637, in margins of Diophantus' *Arithmetica*, was proved successfully by English mathematician, Professor Andrew Wiles at Princeton University, USA, in October, 1994. Professor Andrew Wiles used such mathematical tools as the Kolyvagin-Flach method and the Iwasawa theory. By solving the hypothesis of Taniyama-Shimura (raised by Japanese mathematicians Yutaka Taniyama and Coro Shimura in the 1950's), Wiles proved the famous theorem[1]. The success of Wiles has enhanced the morale of the world's mathematicians, but since then, people have expressed regret about this matter. Because this baffling problem looks so simple, it has challenged the intelligence of innumerable mathematicians and amateurs for 358 years. People can't help asking: why can't the numerous world's outstanding mathematicians overcome even such a simple problem? The proverb "when meet mountain must drill rock, when meet river must erect bridge" seems apposite. Based on research of *n*-dimensional cube and in the light of the proof of Pythagoras theorem, I finished the proof in a short time. Although it has nothing to do with Professor Andrew Wiles' many mathematical tools, and nothing to do with his predecessors' methods, I think the proof is clear at a glance, the tools used are basic, and the result is satisfactory.

## 1 Equivalent Proposition of Fermat's Last Theorem

**Fermat's Last Theorem** *Suppose that $n$ is the positive integers $(n \geq 3)$, then the indefinite equation $x^n + y^n = z^n$ has no any solutions of positive integers.*

Denoting $x$, $y$, $z$ are respectively the lengths of the edges of the distinct three $n$-dimensional cubes, then the above proposition becomes:

**Equivalent Proposition** *Suppose that $n \geq 3$, then the sum of the volumes of two distinct n-dimensional cubes (whose the lengths of the edges are respectively $x$ (positive integer) and $y$ (positive integer)) doesn't equal to the volume of third n-dimensional cube (whose the length of the edge is $z$ (positive integer)).*

This Proposition corresponds to the following figure:

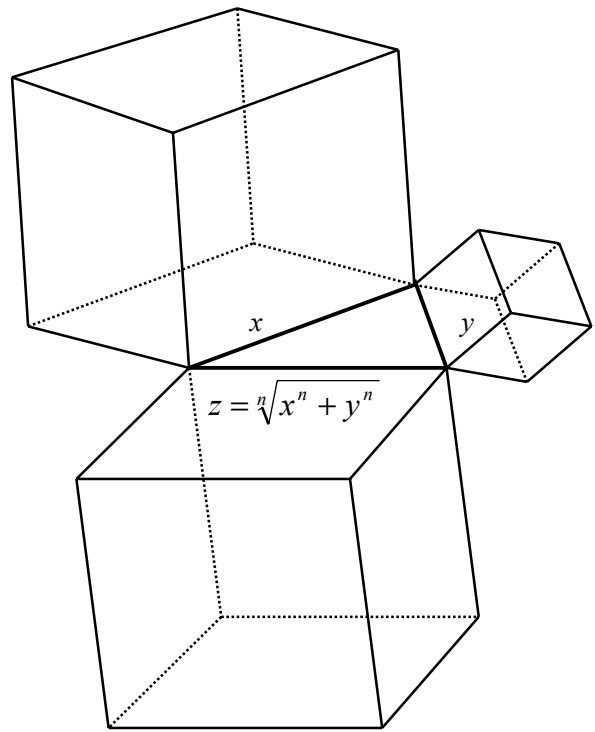

Figure 1

## 2 Using "Embedding" Method to Prove Fermat's Last Theorem

For proving the indefinite equation $x^2 + y^2 = z^2$ has infinite groups of the solutions of positive integers, **Pythagoras** used the following figure:

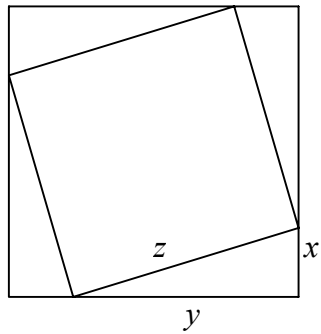

Figure 2

Where, $x$, $y$, $z$ are respectively the lengths of three borders of a rectangular triangle. The geometric meaning of **Pythagoras theorem** is:

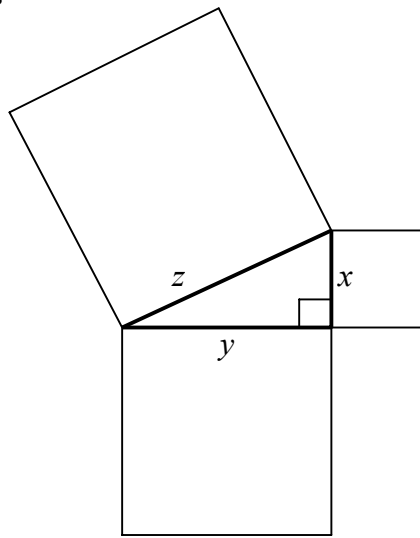

Figure 3

For the indefinite equation $x^2 + y^2 = z^2$, people have found: in following solution

$$\begin{cases} x = u^2 - v^2, \\ y = 2uv, \quad (u > v,\ u,\ v \in Z^+) \\ z = u^2 + v^2 \end{cases} \quad (1)$$

when $u,\ v$ are taken respectively the certain positive integers, the $x,\ y,\ z$ will obtain infinite groups of solutions of positive integers.

Now, let's make a change for the indefinite equation $x^n + y^n = z^n$, that is

$$z = \sqrt[n]{x^n + y^n}.$$

We can know that when $n = 3,\ 4,\ \cdots,\ \sqrt{x^2 + y^2} > \sqrt[n]{x^n + y^n}$, therefore, the angle $\theta$ (which corresponds to $z = \sqrt[n]{x^n + y^n}\ (n = 3,\ 4,\ \cdots)$) less than the rectangular angle $\dfrac{\pi}{2}$ (which corresponds to $z = \sqrt{x^2 + y^2}$), namely $\theta < \dfrac{\pi}{2}$, therefore $\cos\theta > 0$.

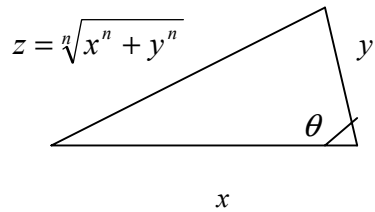

Figure 4

According to the **Cosine Theorem**, we have

$$\left(\sqrt[n]{x^n + y^n}\right)^2 + 2xy\cos\theta = x^2 + y^2. \quad (2)$$

When $n = 2$, because $\theta = \dfrac{\pi}{2}$, $2xy\cos\theta = 0$, therefore, (2) becomes following identity:

$$\left(\sqrt{x^2 + y^2}\right)^2 = x^2 + y^2.$$

When $n \geq 3$, the triangle (Figure 4) takes respectively $x$, $y$, $z = \sqrt[n]{x^n + y^n}$ as the lengths of three borders. In order to make a deep research in it, we must use (2) to change it from a non-rectangular triangle (Figure 4) to a rectangular triangle (Figure 5).

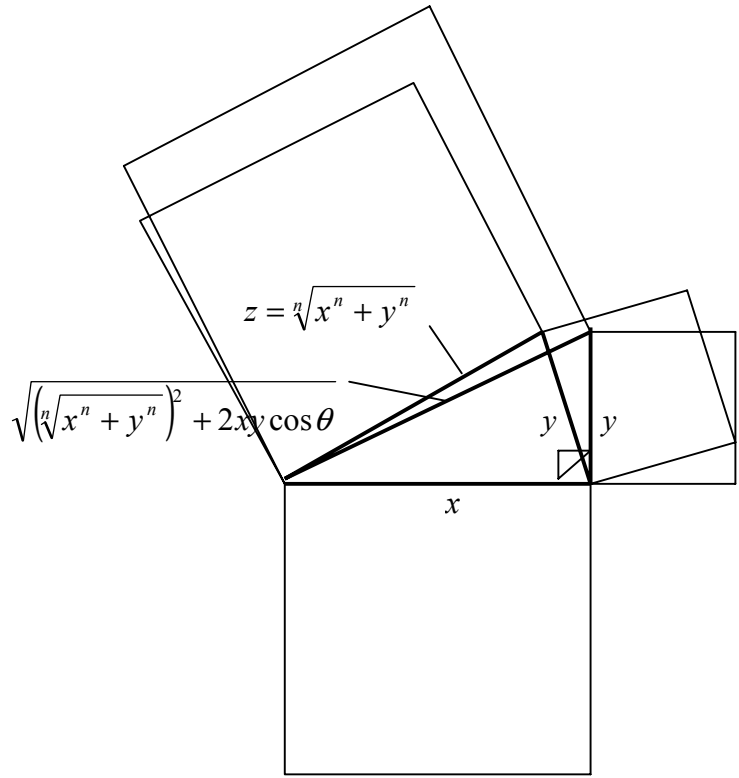

Figure 5

The **key** of **Fermat's Last Theorem** is just in the **triangle**.

For the rectangular triangle ($n = 2$, $\theta = \frac{\pi}{2}$), we have known that the indefinite equation $x^2 + y^2 = z^2$ has infinite groups of solutions of positive integers; but for the non-rectangular triangle ($n \geq 3$, $\theta < \frac{\pi}{2}$), because the triangle doesn't accord with the condition (rectangular triangle) of **Pythagoras theorem**, therefore, according to the **Cosin theorem**, we must add **a line segment** (whose length is $2xy\cos\theta$) to the square number of the original line segment (whose value is $z = \sqrt[n]{x^n + y^n}$). So, we may obtain a rectangular triangle (Figure 5). Obviously, the length of the hypotenuse of the rectangular triangle is $\sqrt{\left(\sqrt[n]{x^n+y^n}\right)^2+2xy\cos\theta}$. According to **Pythagoras theorem**, we have

$$\left(\sqrt{\left(\sqrt[n]{x^n+y^n}\right)^2+2xy\cos\theta}\right)^2 = x^2+y^2.$$

When $n=2$, $\theta=\dfrac{\pi}{2}$, the term $2xy\cos\theta = 0$, the **overall size of the hypotenuse** of the rectangular triangle is just

$$\sqrt{\left(\sqrt{x^2+y^2}\right)^2+2xy\cos\theta} = \sqrt{x^2+y^2}.$$

So, the solution (1) can be embedded into the equation $\left(\sqrt{x^2+y^2}\right)^2 = x^2+y^2 = z^2$. Hence, the indefinite equation $x^2+y^2=z^2$ has infinite groups of the solutions of positive integers.

When $n\geq 3$, $\theta<\dfrac{\pi}{2}$, the term $2xy\cos\theta>0$. Because the **overall size of the hypotenuse** of the rectangular triangle is **invariant**, therefore the term $\sqrt[n]{x^n+y^n}$ has been **reduced**. So, two quantities $x=u^2-v^2$, $y=2uv$ of the solution (1) keep their original positions, although the third quantity $z=u^2+v^2$ of the solution (1) **can be embedded into** the **overall length** $\sqrt{\left(\sqrt[n]{x^n+y^n}\right)^2+2xy\cos\theta}$ **of the hypotenuse** of the rectangular triangle, but **cannot be embedded into** the term $z=\sqrt[n]{x^n+y^n}$. On the other hand, considering the lengths of three borders of the triangle are respectively the lengths of the edges of three *n*-dimensional cubes (whose volumes are respectively $x^n$, $y^n$, $z^n$), therefore, however, they are just the solutions of the indefinite equation $x^n+y^n=z^n$. Hence, we can know when $n\geq 3$, the indefinite equation $x^n+y^n=z^n$ has no any solutions of positive integers.

## 3 Explanation

Maybe, people would like to point out that according to the Cosine theorem, in the equation (2)

$$\left(\sqrt[n]{x^n + y^n}\right)^2 + 2xy\cos\theta = x^2 + y^2,$$

the three quantities $\sqrt[n]{x^n + y^n}$, $x$ and $y$ can be taken positive integers, but we can know that there is a very important relationship exists in them, namely:

$$x^n + y^n = z^n,$$

this is just the problem of we want to prove, therefore, the three quantities $\sqrt[n]{x^n + y^n}$, $x$ and $y$ cannot be taken positive integers in advance.

By researching the *n*-dimensional cube, we found a new way (geometric method), and have proved the world's mathematical baffling problem, Fermat's Last Theorem, we are also known that Fermat's conjecture is true.

1994 (1).

[7] P. Ribenboim, 13 *Lectures on Fermat's Last Theorem* , Springer-Verlag, 1979.

# Appendix

## Dear Mr. Referees,

**Please observe the following photos of the 4-dimensional cube, 7-dimensional cube and the table of orbit radix of the cubes in $E^4 \sim E^7$:**

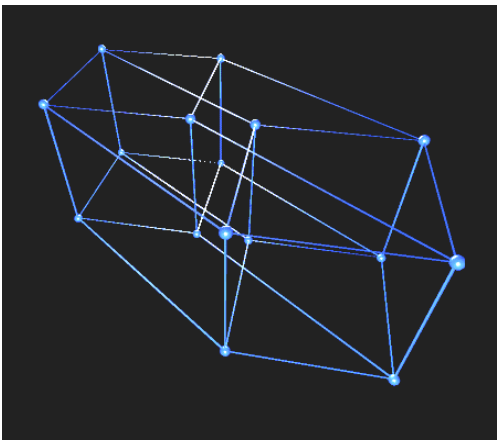
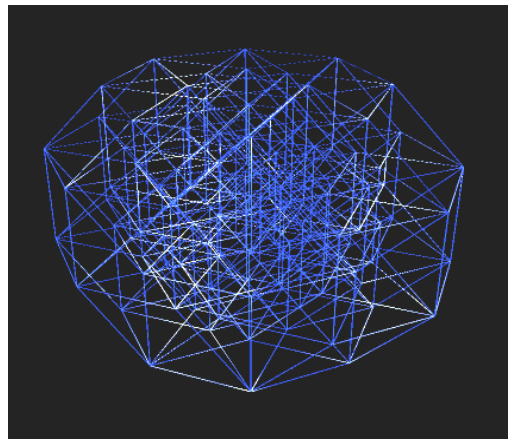

4-dimensional cube                    7-dimensional cube

### Table of Orbit Radix of the Cubes in $E^4 \sim E^7$

| Name     | Vertices | Edges | Squares | Cubes | 4-D Cubes | 5-D Cubes | 6-D Cubes |
|----------|----------|-------|---------|-------|-----------|-----------|-----------|
| 4-D Cube | 16       | 32    | 24      | 8     |           |           |           |
| 5-D Cube | 32       | 80    | 80      | 40    | 10        |           |           |
| 6-D Cube | 64       | 192   | 240     | 160   | 60        | 12        |           |
| 7-D Cube | 128      | 448   | 672     | 560   | 280       | 84        | 14        |

## *Dr. Kaida Shi,*
*Department of Mathematics,*

*Zhejiang Ocean University,*
*CHINA.*